\documentclass[10pt]{article}

\usepackage{graphics}
\usepackage{graphicx}

\usepackage[a4paper, left=35mm,right=35mm,top=34mm,bottom=34mm]{geometry}
\usepackage[utf8]{inputenc}
\usepackage[T1]{fontenc}
\usepackage[english]{babel}

\usepackage{enumerate}
\usepackage{graphicx}
\usepackage{hyperref}
\hypersetup{
    colorlinks=true,
    linkcolor=blue,
    filecolor=magenta,      
    urlcolor=cyan,
}

\usepackage{mathtools,amsthm,amssymb,amsfonts}
\usepackage{algorithm}
\usepackage{algorithmic}
\makeatother
\theoremstyle{plain}

\theoremstyle{definition}

\theoremstyle{remark}

\usepackage{caption} 
\usepackage[tight,footnotesize]{subfigure}

\usepackage{fancyhdr}

\author{
  {\normalsize Abal-Kassim Cheik Ahamed}\thanks{\'Ecole Centrale Paris, France.}
	\and
  {\normalsize Fr\'ed\'eric Magoul\`es}\thanks{\'Ecole Centrale Paris, France
    (correspondence, frederic.magoules@hotmail.com).}
		}	
\title{Stochastic Optimized Schwarz Methods for the Gravity Equations on Graphics Processing Unit}
\date{}

\begin{document}
\maketitle
\thispagestyle{fancy}

\begin{abstract}
\noindent Low order, sequential or non-massively parallel finite elements are generaly used for three-dimensional gravity modelling.
In this paper, in order to obtain better gravity anomaly solutions in heterogeneous media, we solve the gravimetry problem using massively parallel high order finite elements on hybrid multi-CPU/GPU clusters.
Parallel algorithms well suited for such hybrid architectures have to be designed.
A new stochastic-based optimization procedure for the optimized Schwarz method is here presented, implemented and tuned to graphical cards processors units.
Numerical experiments performed on a reallistic test case, demonstrates the robustness and efficiency of the proposed method and of its implementation on massive multi-CPU/GPU architectures.

\end{abstract}


\section{Introduction}\label{magoulesf_contrib_3:sec:1}

By giving another way to see beneath the Earth, gravimetry refines geophysical exploration.
In this paper, we evaluate the gravimetry field in the Chicxulub crater area located in between the Yucatan region and the Gulf of Mexico which shows strong gravimetry and magnetic anomalies.
High order finite elements analysis is considered with input data arising from real measurements.
The linear system is then solved with a domain decomposition method, namely the optimized Schwarz method.
The principle of this method is to decompose the computational domain into smaller subdomains and to solve the equations on each subdomain.
Each subdomain could easily be allocated to one single processor (i.e. the CPU), each iteration of the optimized Schwarz method involving the solution of the equations on each subdomain (on the GPU).
Unfortunately, to obtain high speed-up, several tunings and adaptations of the algorithm should be carrefully performed, such as data transfers between CPU and GPU, and data structures, as described in ~\cite{magoulesf_contrib_3:CheikAhamedMagoules:2012,magoulesf_contrib_3:CheikAhamedMagoules:2012:a}.

The plan of the paper is the following.
In Section~\ref{magoulesf_contrib_3:sec:2}, we present the gravimetry equations.
In Section~\ref{magoulesf_contrib_3:sec:3}, we introduce the optimized Schwarz method, followed in Section~\ref{magoulesf_contrib_3:sec:4} by a new idea of using a stochastics-based algorithm to determine the optimized transmission conditions.
An overview to the GPU programming model and hardware configuration suite is given in Section~\ref{magoulesf_contrib_3:sec:5} for readers not familiar with GPU programming.
Section~\ref{magoulesf_contrib_3:sec:6} shows numerical experiments which clearly confirm the robustness, competitiveness and efficiency of the proposed method on GPU clusters for solving the gravimetry equations.

\section{Gravimetry equations}\label{magoulesf_contrib_3:sec:2}

The gravity force is the resultant of the gravitational force and the centrifugal force.
The gravitational potential of a spherical density distribution is equal to
$\Phi(r) = Gm / r$, with $m$ the mass of the object, $r$ the distance to the object and $G$ the universal gravity constant equal to $G = 6.672\times 10^{-11} m^{3}kg^{-1}s^{-2}$.
The gravitational potential at a given position $x$ initiated by an arbitrary density distribution $\rho$ is given by $\Phi(x) = G\int (\rho(x') / ||x-x'||) dx'$
where $x'$ represents one point position within the density distribution.
In this paper, we consider only regional scale of the gravimetry equations therefore we do not take into account the effects related to the centrifugal force.
The gravitational potential $\Phi$ of a density anomaly distribution $\delta\rho$ is thus given as the solution of the Poisson equation $\Delta \Phi = -4\pi G \delta \rho$.

\section{Optimised Schwarz method}\label{magoulesf_contrib_3:sec:3}

The classical Schwarz algorithm was invented more than a century ago~\cite{magoulesf_contrib_3:HASchwarz} to prove existence and uniqueness of solutions to Laplace's equation on irregular domains.
Schwarz decomposed the irregular domain into overlapping regular ones and formulated an iteration which used only solutions on regular domains and which converged to a unique solution on the irregular domain.
The convergence speed of the classical Schwarz algorithm is proportional to the size of the overlap between the subdomains.
A variant of this algorithm can be formulated with non-overlapping subdomains and the transmission conditions should be changed from Dirichlet to Robin~\cite{magoulesf_contrib_3:Despres:1992:DHM}.
These absorbing boundary transmission conditions defined on the interface between the non-overlapping subdomains, are the key ingredients to obtain a fast convergence of the iterative Schwarz algorithm~\cite{magoulesf_contrib_3:Despres:1993:DDM,magoulesf_contrib_3:Benamou:1997:DDM}.
Optimal transmission conditions can be derived but consists of non local operators and thus are not easy to implement in a parallel computational environment.
One alternative is to approximate these operators with partial differential operators.
This paper investigates an approximation based on a new stochastics optimization procedure.

For the sake of clarity, the gravimetry equations are considered in the domain $\Omega$ with homogeneous Dirichlet condition.
The domain is decomposed into two non-overlapping subdomains $\Omega^{(1)}$ and $\Omega^{(2)}$ with an interface $\Gamma$.
The Schwarz algorithm can be written as:
\begin{eqnarray*}
-\Delta\Phi^{(1)}_{n+1} & = & f,
\quad \mbox{in} \quad \Omega^{(1)} \\
\left(\partial_\nu \Phi^{(1)}_{n+1}+\mathcal{A}^{(1)} \Phi^{(1)}_{n+1}\right) & = &
\left(\partial_\nu \Phi^{(2)}_{n}+\mathcal{A}^{(1)} \Phi^{(2)}_{n}\right),
\quad \mbox{on} \quad \Gamma \\
-\Delta\Phi^{(2)}_{n+1} & = & f,
\quad \mbox{in} \quad \Omega^{(2)} \\
\left(\partial_\nu \Phi^{(2)}_{n+1}-\mathcal{A}^{(2)} \Phi^{(2)}_{n+1}\right) & = &
\left(\partial_\nu \Phi^{(1)}_{n}-\mathcal{A}^{(2)} \Phi^{(1)}_{n}\right),
\quad \mbox{on} \quad \Gamma
\end{eqnarray*}
with $n$ the iteration number, and $\nu$ the unit normal vector defined on $\Gamma$.
The operators $\mathcal{A}^{(1)}$ and $\mathcal{A}^{(2)}$ are to be determined for best performance of the algorithm.
Considering the case $\Omega=\mathbb{R}^2$, $f=0$, and applying a Fourier transform, similar calculations as in~\cite{magoulesf_contrib_3:Gander:2002:OSM} lead to the expression of the Fourier convergence rate, involving the quantities  $\Lambda^{(1)}$ and $\Lambda^{(2)}$, which are the Fourier transforms of $\mathcal{A}^{(1)}$ and $\mathcal{A}^{(2)}$ operators.
In this case, the choice $\Lambda^{(1)}:=|k|$, and $\Lambda^{(2)}:=|k|$ is optimal, since with this choice the algorithm converges in two iterations for two subdomains.
Different techniques to approximate these non local operators with partial differential operators have been analyzed in recent years~\cite{magoulesf_contrib_3:Despres:1993:DDM,magoulesf_contrib_3:Chevalier:1998:SMO,magoulesf_contrib_3:Gander:2002:OSM}.
These techniques consist to define partial differential operators involving a tangential derivative on the interface such as: ${\cal A}^{(s)} := p^{(s)} + q^{(s)} \partial^2_{\tau^2}$, with $s$ the subdomain number, $p^{(s)}$, $q^{(s)}$ two coefficients, and $\tau$ the unit tangential vector.
The first results presented in~\cite{magoulesf_contrib_3:Despres:1993:DDM,magoulesf_contrib_3:Benamou:1997:DDM} use a zero order Taylor expansion of the non local operators to find $p^{(s)}$ and $q^{(s)}$.
In~\cite{magoulesf_contrib_3:Japhet:2000:OO2} for convection diffusion equations, in~\cite{magoulesf_contrib_3:Chevalier:1998:SMO} for Maxwell equation, in~\cite{magoulesf_contrib_3:Gander:2002:OSM,magoulesf_contrib_3:Magoulestopping2} for the Helmholtz equation,
and in~\cite{magoulesf_contrib_3:Madaymagoules,magoulesf_contrib_3:Madaymagoules2} for heterogeneous media, a minimization procedure has been used.
The function to minimize, i.e., the cost function, is the maximum of the Fourier convergence rate for the frequency ranges considered, and the approach consists to determine the free parameters $p^{(s)}$ and $q^{(s)}$ through an optimization problem.
Despite, analytic expressions of $p^{(s)}$ and $q^{(s)}$ can be determined for some specific problems, finding quasi-optimal coefficients numerically is also a good alternative~\cite{magoulesf_contrib_3:Magoulestopping}.
Furthemore, since the evaluation of the cost function is quite fast and the dimension of the search space reasonable, a more robust minimization procedure could be considered, in the next section. Extension to non-regular geometry can be performed as described in reference~\cite{magoulesf_contrib_3:MagoulesPutanowicz:2005:DDM}.

\section{Stochastic-based optimised transmission conditions}\label{magoulesf_contrib_3:sec:4}

The stochastic minimization technique we propose to use now, explores the whole space of solutions and finds absolute minima; this technique is based on the Covariance Matrix Adaptation Evolution Strategy (CMA-ES).
This algorithm is very robust~\cite{magoulesf_contrib_3:AugerH:2012}, has good global search ability and does not need to compute the derivatives of the cost function.
This algorithm only needs an initial search zone and a population size, even if the solution can be found outside of the initial search zone.
The population size parameter is a trade-off between speed and global search.
Meaning that, smaller populations lead to faster execution of the algorithm but have more chance to find a local minimum, and, larger sizes allow to avoid local minima better but need more cost function evaluations.
For our purpose, a population size of 25 has been large enough to find the global minimum in a few second or less.

The main idea of the algorithm is to find the minimum of the cost function by iteratively refining a search distribution.
The distribution is described as a general multivariate normal distribution $d(m,C)$.
Initially, the distribution is given by the user.
Then, at each iteration, $\lambda$ samples are randomly chosen in this distribution and the evaluation of the cost function at those points is used to compute a new distribution.
When the variance of the distribution is small enough, the center of the distribution $m$ is taken as solution.
After evaluating the cost function for a new population, the samples are sorted by cost and only the $\mu$ best are kept.
The new distribution center is computed with a weighted mean (usually, more weight is put on the best samples).
The step size $\sigma$ is a factor used to scale the standard deviation of the distribution, i.e., the variance of the search distribution is proportional to the square of the step size.
The step size determines the ``size'' of the distribution.
The covariance matrix $C$ determines the ``shape'' of the distribution, i.e., it determines the principal directions of the distribution and their relative scaling.
Adapting (or updating) the covariance matrix is the most complex part of the algorithm.
While this could be done using only the current population, it would be unreliable especially with a small population size; thus the population of the previous iteration should also been taken into account.

\section{Overview of GPU programming model}\label{magoulesf_contrib_3:sec:5}

Parallel computation was generally carried out on Central Processing Unit (CPU) cluster until the apparition in the early 2000s of the Graphics Processing Unit (GPU) that facing the migration of the era of GPU computing.
The peak performance of CPUs and GPUs is significanlty different, due to the inherently different architectures between these processors.
The first idea behind the architecture of GPU is to have many small floating points processors exploiting large amount of data in parallel.
This is achieved through a memory hierarchy that allows each processor to optimally access the requiered data.
The gains of GPU computing is significantly higher for large size problem compared to CPU, due to the difference between these two architectures.
GPU computing requires using graphics programming languages such as NVIDIA CUDA, or OpenCL.
\emph{Compute Unified Device Architecture} (CUDA)~\cite{magoulesf_contrib_3:CUDA4.0} programming model is an extension of the C language and has been used in this paper.

A specific characteristic of GPU compared to CPU is the feature of memory used.
Indeed, a CPU is constantly accessing the RAM, therefore it has a low latency at the detriment of its raw throughput.
CUDA devices have four main types of memory: (i) \emph{Global memory} is the memory that ensures the interaction with the host (CPU), and is not only large in size and off-chip, but also available to all threads (also known as compute units), and is the slowest; (ii) \emph{Constant memory} is read only from the device, is generally cached for fast access, and provides interaction with the host;
(iii) \emph{Shared memory} is much faster than global memory and is accessible by any thread of the block from which it was created;
(iv) \emph{Local memory} is specific to each compute unit and cannot be used to communicate between compute units.

\begin{figure}  
\centerline{\subfigure{\includegraphics[scale=0.26]{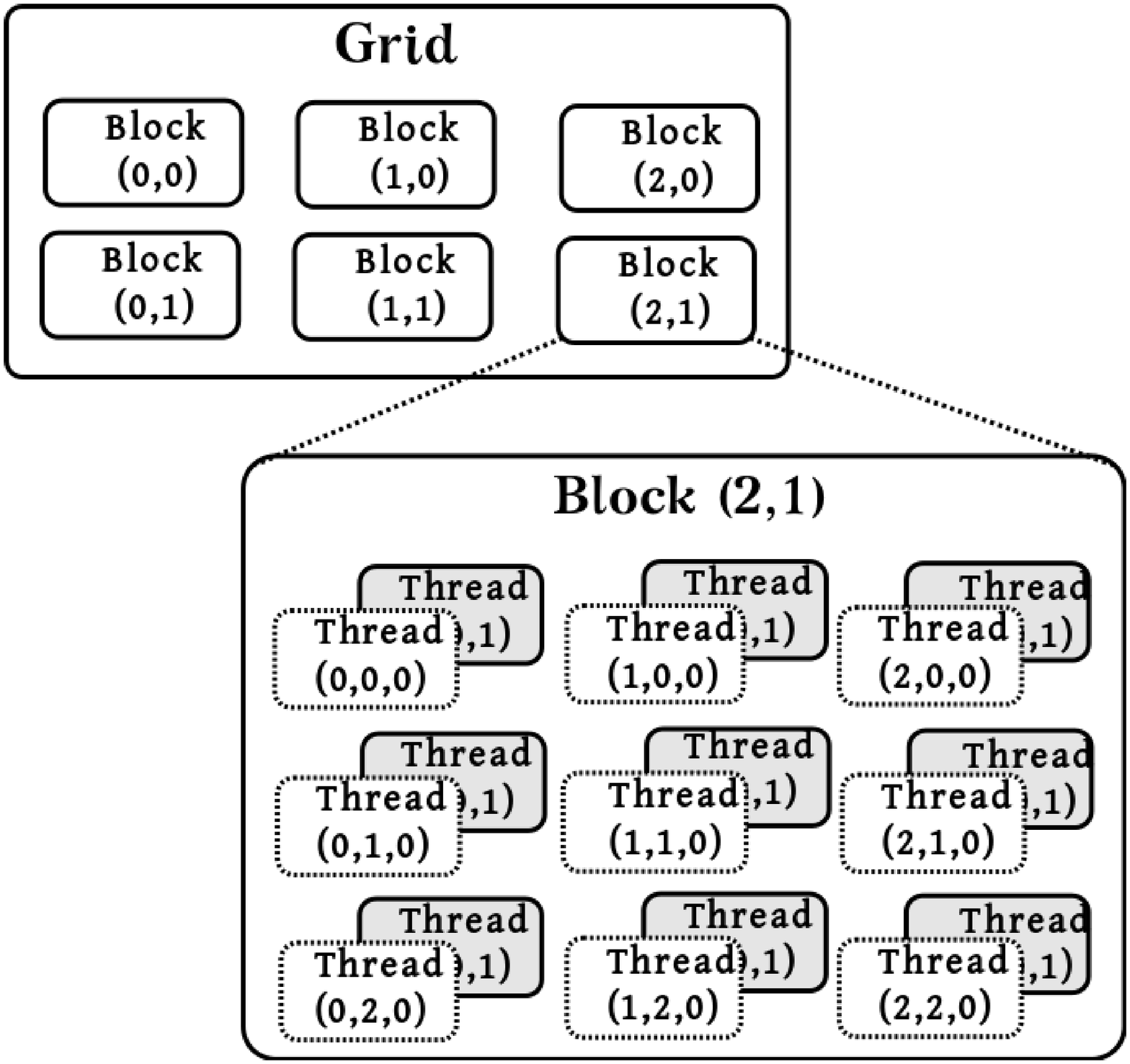}}
\hfil  
\subfigure{\includegraphics[scale=0.26]{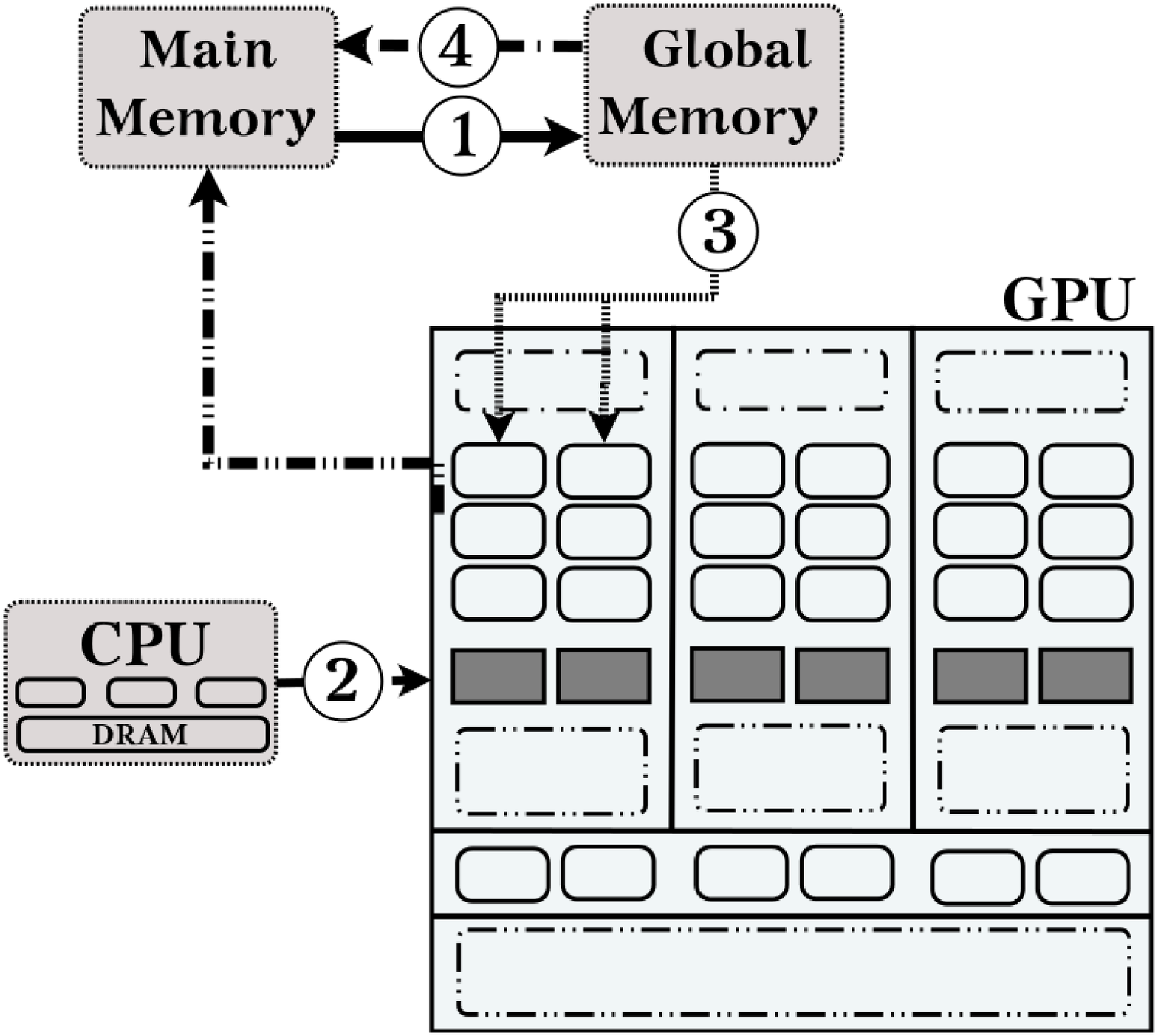}}} 
 \caption{Gridification of a GPU. Tthread, block, grid (left); GPU computing processing (right)}
 \label{magoulesf_contrib_3:fig:img:flow_cuda}
\end{figure}

Threads are grouped into blocks and executed in parallel simultaneously, see Figure~\ref{magoulesf_contrib_3:fig:img:flow_cuda}.
A GPU is associated with a \emph{grid}, i.e., all running or waiting blocks in the running queue and a kernel that will run on many cores.
An ALU is associated with the thread which is currently processing.
Threading is not an automated procedure. The developer chooses for each kernel the distribution of the threads, which are organized (\emph{gridification} process) as follows: (i) threads are grouped into blocks; (ii) each block has three dimensions to classify threads; (iii) blocks are grouped together in a grid of two dimensions.
The threads are then distributed to these levels and become easily identifiable by their positions in the grid according to the block they belongs to and their spatial dimensions.
The kernel function must be called also providing at least two special parameters: the dimension of the block, \emph{nBlocks}, and the number of threads per block, \emph{nThreadsPerBlock}.
Figure~\ref{magoulesf_contrib_3:fig:img:flow_cuda} presents the CUDA processing flow. Data are first copied from the main memory to the GPU memory, (1).
Then the host (CPU) instructs the device (GPU) to carry out calculations, (2).
The kernel is then executed by all threads in parallel on the device, (3).
Finally, the device results are copied back (from GPU memory) to the host (main memory), (4).
\\
To cope with this difficulty the implementation proposed in this paper uses some advanced tuning techniques developped by the authors, but the details are outside the scope of this paper, and the reader is refeered to~\cite{magoulesf_contrib_3:CheikAhamedMagoules:2012,magoulesf_contrib_3:CheikAhamedMagoules:2012:a} for the computer science aspects of this tuning.

\section{Numerical analysis}\label{magoulesf_contrib_3:sec:6}

In this section, we report the experiments performed to evaluate the speed-up of our implementation.
The Chicxulub impact crater, formed about 65 million years ago, is now widely accepted as the main footprint of the global mass extinction event that marked the Cretaceous/Paleogene boundary.
Because of its relevance, in the last two decades, this impact structure has been used as a natural laboratory to investigate impact cratering formation processes and to indirectly infer global effects of large-scale impacts.
The crater is buried under 1 km of carbonate sediments in the Yucatan platform. The crater is about 200 km in rim diameter, half on-land and half off-shore with geometric center at Chicxulub Puerto.
The internal structure of the Chicxulub crater has been imaged by using several geophysical data sets from land, marine and aerial measurements.

In this paper we perform a finite element analysis of the gravimetry equation using the characteristics of the region provided by the measure.
The domain consists of an area of $250 \times 250 \times 15$ in each spatial direction, and is discretized with high order finite element with a total of 19~933~056 degrees of freedom. 
This mesh is partionned in the x-direction. 
Each subdomain is allocated to one single processor (i.e. the CPU), each iteration of the optimized Schwarz method involving the solution of the equations inside each subdomain is allocated to one single accelerator (i.e the GPU). 
We compare the computational time of the optimised Schwarz method using one subdomain per CPU with the optimised Schwarz method using one subdomain per CPU and a GPU accelerator.
For this particular model we performed calculations using our CUDA implementation of the Schwarz method with stochastics-based optimization procedure.
The workstation used for all the experiments consists of 1~596 servers Bull Novascale R422Intel Nehalem-based nodes.
Each node is composed of 2 processors Intel Xeon 5570 quad-cores (2.93 GHz) and 24 GB of memory (3Go per cores).
96 CPU servers are interconnected with 48 compute Tesla S1070 servers NVIDIA (4 Tesla cards with 4GB of memory by server) and 960 processing units are available for each server.\\

For the subdomain problems, the diagonal preconditioner conjugate gradient (PCG) is used and the coefficient matrices are stored in CSR format. We fix a residual tolerance threshold of $\epsilon = 10^{-10}$ for PCG. \emph{Alinea}~\cite{magoulesf_contrib_3:CheikAhamedMagoules:2012,magoulesf_contrib_3:CheikAhamedMagoules:2012:a}, our research group library, implemented in \emph{C++}, which offers CPU and GPU solvers, is used for solving linear algebra system.
In this paper, the GPU is used to accelerate the solution of PCG algorithm.
PCG algorithm required the computation of addition of vectors (\emph{Daxpy}), dot product and sparse matrix-vector multiplication. In GPU-implementation considered (Alinea library), the distribution of threads (\emph{gridication}, differs with these operations.
The gridification of \emph{Daxpy}, dot product and sparse matrix-vector product correspond respectively to ($nBlocks=\frac{ numb\_rows + numb\_th\_block - 1 }{numb\_th\_block}$, $nThreadsPerBlock=256$), ($nBlocks=\frac{ numb\_rows + numb\_th\_block - 1 }{numb\_th\_block}$, $nThreadsPerBlock=128$) and ($nBlocks=\frac{ (numb\_rows \times n\_th\_warp) + numb\_th\_block - 1 }{numb\_th\_block}$, $nThreadsPerBlock=256$), where $numb\_rows$, $n\_th\_warp$ and $numb\_th\_block$ represent respectively  the number of rows of the matrix, the number of threads per warp and the thread block size. We fix for all operations 8 threads per warp. GPU is used only for solving subdomain problems in each iteration.
GPU experiment workstation Tesla S1070 has 4 GPUs of 240 cores. The number of computing units depends both on the size of the subdomain problem and the gridification that use $256$ threads per threads and $8$ threads per warp as introduced in~\cite{magoulesf_contrib_3:CheikAhamedMagoules:2012,magoulesf_contrib_3:CheikAhamedMagoules:2012:a}.\\

In our experiments, the CMA-ES algorithm considers as the cost function the Fourier convergence rate of the optimised Schwarz method. We consider for the CMA-ES the following stopping criteria of : a maximum of number iterations equal to $7200$ and a residu threshold equal to $5\times10^{-11}$. 
\begin{figure}
\centering
\scalebox{0.5}{\includegraphics{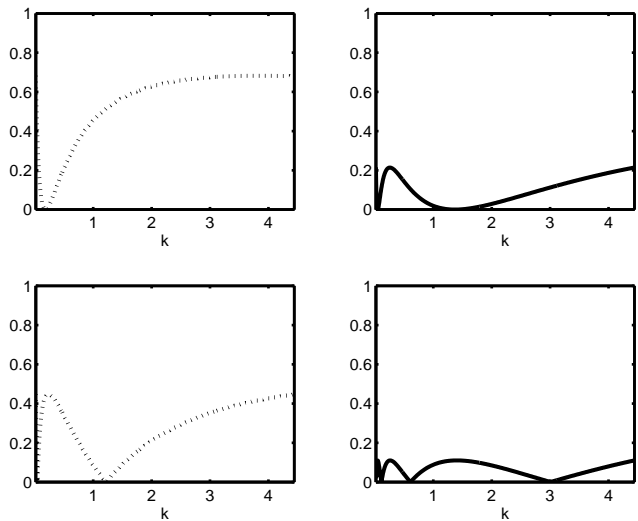}}
\caption{Fourier convergence rate of the Schwarz algorithm}
\label{magoulesf_contrib_3:fig:fourier}
\end{figure}
Fig.~\ref{magoulesf_contrib_3:fig:fourier} represents the convergence rate of the Schwarz algorithm in the Fourier space, respectively for the symmetric zeroth order (top-left), unsymmetric zeroth order (bottom-left), the symmetric second order (top-right) and unsymmetric second order (bottom-right) transmission conditions obtained from the CMA-ES algorithm.
The Fourier convergence rate of the Schwarz method with one side (respectilely two sides) transmission conditions obtained from the CMA-ES algorithm is presented in Fig.~\ref{magoulesf_contrib_3:fig:fourier} and Table~\ref{magoulesf_contrib_3:tab:ddm:opt_configuration}.
\begin{table}
\centering
{
\begin{tabular}{|l|c|c|c|c|c|}
\hline 
 & ${p^{(1)}}$ & ${q^{(1)}}$ & ${p^{(2)}}$ & ${q^{(2)}}$ & ${\rho_{max}}$ \\
\hline %
{oo0\_symmetric}   & 0.1826 & 0 & 0.1826 & 0 & 0.6823 \\
{oo0\_unsymmetric} & 1.2193 & 0 & 0.0469 & 0 & 0.4464 \\
{oo2\_symmetric}   & 0.0471 & 0.7050 & 0.0471 & 0.7050 & 0.2143 \\
{oo2\_unsymmetric} & 0.1081 & 0.3205 & 0.0231 & 1.5786 & 0.1101 \\
\hline 
\end{tabular}
}
\caption{Optimized coefficients obtained from \emph{CMA-ES} algorithm}
\label{magoulesf_contrib_3:tab:ddm:opt_configuration}
\end{table}
\begin{table}
\centering
{
\begin{tabular}{|c|c|c|c|c|}
\hline 
\textbf{\#subdomains} & \textbf{\#iterations}
& \textbf{cpu time (sec)} & \textbf{ gpu time (sec)} &  \textbf{SpeedUp}\\ 
\hline \hline %
32  & 41 & 11~240 & 1~600 & \textbf{7.03}\\
64  & 45 &  5~360 & 860 & \textbf{6.23}\\
128 & 92 & 6~535 & 960 & \textbf{6.81}\\
\hline 
\end{tabular}
}
\caption{Comparison of the implementation of our method on CPU and GPU}
\label{magoulesf_contrib_3:tab:ddm:synch_cg_res_case1oo2_procs_2_nodes}
\end{table}
The distribution of processors is computed as follows: $\text{number of processors} = 2\times \text{number of nodes}$, where $2$ corresponds to the number of GPU per node as available on our workstation. As a consequence, only two processors will share the bandwidth, which strongly improve the communications, especially the inter-subdomain communications.
Table~\ref{magoulesf_contrib_3:tab:ddm:synch_cg_res_case1oo2_procs_2_nodes} presents the results done with double precision with a residu threshold, \emph{i.e.} stopping criterion equal to $10^{-6}$, for several number of subdomains (one subdomain per processor).

\section{Conclusion}\label{magoulesf_contrib_3:sec:7}

In this paper, we have presented a stochastic-based optimized Schwarz method for the gravimetry equation on GPU Clusters.
The effectiveness and robustness of our method are evaluated by numerical experiments performed on a cluster composed of 1~596 servers Bull Novascale R422Intel Nehalem-based nodes where 96 CPU servers are interconnected with 48 compute Tesla S1070 servers NVIDIA (4 Tesla cards with 4GB of memory by server).
The presented results range from 32 up-to 128 subdomains show the interest of the use of GPU technologies for solving large size problems, and outline the robustness, performance and efficiency of our Schwarz domain decomposition method with stochastic-based optimized conditions.

\section*{Acknowledgments}
The authors acknowledge partial financial support from the OpenGPU project (2010-2012), and GENCI (Grand Equipement National de Calcul Intensif) for the computer time used during this long-term trend.

\bibliography{bib/magoulesf_contrib_3}
\bibliographystyle{abbrv}

\end{document}